\documentstyle{amsppt}
\overfullrule=0pt
\font\sc=cmcsc10
\magnification=\magstep1
 \pagewidth{6.0 true in}
\pageheight{9.0 true in}

\hoffset 0in 
\voffset -0.10in

\parindent 0pt
\parskip 14pt
\topmatter
\title Obstructions to generic embeddings
\endtitle
\author Judith Brinkschulte \\
C. Denson Hill \\
Mauro Nacinovich
\endauthor
\endtopmatter
In Grauert's paper [G] it is noted that finite dimensionality of cohomology groups sometimes implies vanishing of these cohomology groups. Later on Laufer formulated a zero-or infinity law for the cohomology groups of domains in Stein manifolds. In this paper we generalize the Laufer's Theorem in [L] and its
version for small domains of $CR$ manifolds, proved in [Br],
by considering Whitney cohomology on locally closed subsets and
cohomology with supports for currents. 
With this approach we obtain a global result for $CR$ manifolds
generically embedded in a Stein manifold. Namely a necessary condition for global embedding
into a Stein manifold is that the $\bar\partial_M$-cohomology groups must be
either zero or infinite dimensional.
\bigskip
\noindent
{\bf \S 1 An abstract Laufer Theorem}
\smallskip
Let $X$ be a Stein manifold of complex dimension $N$. Let $F$ be
a locally closed subset of $X$. This means that $F$ is a closed
subset of an open submanifold $Y$ of $X$. We denote by $\frak W_F$
the space of {\it Whitney functions} on $F$. With $\frak F(Y;F)$ denoting
the subspace of the space
$\Cal E(Y)$ of (complex valued) smooth functions on $Y$ that vanish
of infinite order at each point of $F$, the space $\frak W_F$
is defined by the exact sequence:
$$\CD
0@>>>\frak F(Y;F)@>>>\Cal E(Y)@>>>\frak W_F@>>>0\, .
\endCD$$
Note that the space $\frak W_F$ can be intrinsically defined in terms
of jets and turns out to be independent of the choice of the open
neighborhood $Y$ of $F$ in $X$. We also consider the space 
$\frak W^{\roman{comp}}_F$ of
{\it Whitney functions with compact support in $F$}, which can be
defined by the exact sequence:
$$\CD
0@>>>\frak F(Y;F)\cap\Cal D(Y)@>>>\Cal D(Y)@>>>
\frak W_F^{\roman{comp}}@>>>0\, ,
\endCD$$
where $\Cal D(Y)$ is the standard notation for the space of $f\in\Cal E(Y)$
having compact support in $Y$.
\par
Likewise we shall consider the spaces $\Cal D'_F$ (resp. $\Cal E'_F$)
of distributions in $Y$ with support (resp. compact support) contained in $F$.
\smallskip
The Dolbeault complexes on $Y$ define, by passing to sub-complexes
and quotients, $\bar\partial$-complexes on {\it Whitney forms} with
closed (or compact) supports in $F$ and on currents with closed (or compact)
supports contained in $F$. We denote by 
$H^{p,q}_{\bar\partial}(\frak W_F)$, 
$H^{p,q}_{\bar\partial}(\frak W^{\roman{comp}}_F)$,
$H^{p,q}_{\bar\partial}(\Cal D'_F)$,
$H^{p,q}_{\bar\partial}(\Cal E'_F)$
the corresponding cohomology groups
(see, for more details, [N1], [N2], [NV]). 
More generally, if $\Phi$ is a paracompactifying family of supports in $Y$ (see [B]),
we can consider the cohomology groups
$H^{p,q}_{\bar\partial}(\frak W_F^{\Phi})$  
for Whitney forms on $F$ with supports in $\Phi$, which 
are quotients of smooth forms in $Y$ with supports in $\Phi$, and 
$H^{p,q}_{\bar\partial}(\Cal D_F^{'\Phi})$ for currents with
supports in the intersection of $F$ and closed sets of $\Phi$.
\par
Note that when $F=Y$ is open, these
are the usual Dolbeault cohomology groups.
\smallskip
\proclaim{Theorem 1.1} Let $F$ be a locally closed subset of a
Stein manifold $X$. Let $Y$ be an open neighborhood of $F$ in $X$,
with $\overline F\cap Y=F$,
and $\Phi$ a paracompactifying family of supports in $Y$. Then,
for all $0\leq p,q\leq N$, each of the cohomology
groups 
$H^{p,q}_{\bar\partial}(\frak W_F)$, 
$H^{p,q}_{\bar\partial}(\frak W^{\roman{comp}}_F)$,
$H^{p,q}_{\bar\partial}(\Cal D'_F)$,
$H^{p,q}_{\bar\partial}(\Cal E'_F)$,
$H^{p,q}_{\bar\partial}(\frak W_F^{\Phi})$,
$H^{p,q}_{\bar\partial}(\Cal D_F^{'\Phi})$,
is either $0$ or infinite dimensional.\endproclaim
\demo{Proof}
The proof follows the argument in [L]. \par
Denote by $\bold H$ one of the
groups $H^{p,q}_{\bar\partial}(\frak W_F)$, 
$H^{p,q}_{\bar\partial}(\frak W^{\roman{comp}}_F)$,
$H^{p,q}_{\bar\partial}(\Cal D'_F)$,
$H^{p,q}_{\bar\partial}(\Cal E'_F)$,
$H^{p,q}_{\bar\partial}(\frak W_F^{\Phi})$,
$H^{p,q}_{\bar\partial}(\Cal D_F^{'\Phi})$,
and assume that $\bold H$ is finite dimensional. Our goal is to show that
$\bold H=\{0\}$. \par
The multiplication of a Whitney form or of a current by a function
$f\in\Cal O(X)$ is a linear map, commuting with $\bar\partial$, and
preserving supports. Thus, by passing to the quotient, we obtain on
$\bold H$ the structure of an $\Cal O(X)$-module.\par
Let $\Cal I=\{f\in\Cal O(X)\, | \, f\bold H=\{0\}\}$ be the ideal
in the ring $\Cal O(X)$ of functions that annihilate $\bold H$.
We want to show that $1\in\Cal I$.
\par
Fix an embedding $X\hookrightarrow\Bbb C^{2N+1}$ of $X$ into a
Euclidean space. The coordinates $z_1,\hdots,z_{2N+1}$ on
$\Bbb C^{2N+1}$ define functions 
$z_1^*,\hdots,z_{2N+1}^*$ in $\Cal O(X)$. Fix $j\in \lbrace 1,\ldots ,2N+1\rbrace$, and let $\lbrack u_1\rbrack,\ldots ,\lbrack u_m\rbrack$ be a basis of $\bold H$. Then the finite dimensionality of $\bold H$ implies that there exist nontrivial polynomials $P_{\lbrack u_i\rbrack}$ such that $P_{\lbrack u_i\rbrack}(z_j^* )\lbrack u_i\rbrack =\lbrace 0\rbrace$, $i= 1,\ldots ,m$. Consider $P=P_{\lbrack u_1\rbrack}\cdot \ldots\cdot P_{\lbrack u_m\rbrack}$. Then one has
$P(z_j^* )\bold H = \lbrace 0\rbrace$. Thus for each $j=1,\hdots,2N+1$
there is a polynomial $P_j(z_j)\in\Bbb C[z_j]\setminus\{0\}$ of minimal
degree such that $P_j(z_j^*)\bold H=\{0\}$. 
This shows on the one hand that $\Cal I\neq\{0\}$ and on the other hand
that the set $V$ of common zeros in $X$ of the functions in $\Cal I$
is finite, being contained in the inverse image by the embedding 
$X\hookrightarrow\Bbb C^{2N+1}$ of the finite set 
$\{z\in\Bbb C^{2N+1}\, | \, P_j(z_j)=0\;\text{for}\; j=1,\hdots,2N+1\}$.
By the Nullstellensatz $1\in\Cal I$ if and only if $V=\emptyset$.
To show this, we prove first the following:\enddemo
\proclaim{Lemma 1.2} Let $f\in\Cal I$ and let $\lambda$ be a
holomorphic vector field on $X$. Then $\lambda(f)\in\Cal I$.\endproclaim
\demo{Proof}
We recall the formula for the Lie derivative 
$L_\lambda$ of an exterior differential
form $\alpha$:
$$L_{\lambda}\alpha=d(\lambda\rfloor\alpha)+\lambda\rfloor d\alpha\, .$$
If $f$ is a smooth function, then
$$L_{\lambda}(f\cdot\alpha)=(\lambda(f))\cdot\alpha+f\cdot L_\lambda\alpha\, .$$
Thus we obtain:
$$\left(\lambda(f)\right)\cdot\alpha=
d(\lambda\rfloor(f\cdot\alpha))+\lambda\rfloor d(f\cdot\alpha)-
f\cdot d(\lambda\rfloor\alpha)-f\cdot\lambda\rfloor d\alpha
\, .\tag \text{$*$} $$
Assume now that $\lambda$ is a holomorphic vector field, that $f$ is
a holomorphic function on $X$ and $\alpha$ is a form 
(or a current) of bidegree
$(p,q)$ in $Y$. Since the left hand side of ($*$) is then of bidegree
$(p,q)$, keeping on the right hand side only the summands which are
homogeneous of bidegree $(p,q)$ we obtain:
$$\left(\lambda(f)\right)\cdot\alpha=
\partial(\lambda\rfloor(f\cdot\alpha))+\lambda\rfloor \partial(f\cdot\alpha)-
f\cdot \partial(\lambda\rfloor\alpha)-f\cdot\lambda\rfloor \partial\alpha
\, .\tag \text{$**$} $$
Assume now that $\alpha$ is $\bar\partial$-closed. Then each term on
the right hand side is $\bar\partial$-closed. By considering subspaces
and quotients, we note that ($**$) is valid 
and that the summands on the right hand side are
$\bar\partial$-closed also when $\alpha$ is
a Whitney form on $F$ or a current with support in $F$.
\par
Assume now that $f\in\Cal I$ and that $\alpha$ is the representative
of an element of $\bold H$. Then the last two summands are cohomologous
to zero because $f\in\Cal I$. Moreover, $f\cdot\alpha=\bar\partial\beta$,
again because $f\in\Cal I$. Thus we have:
$$\matrix\format\l\\
\partial(\lambda\rfloor(f\cdot\alpha))=\partial(\lambda\rfloor\bar\partial\beta)
=\bar\partial\left(\partial(\lambda\rfloor\beta)\right)\, ,\\
\lambda\rfloor\partial(f\cdot\alpha)=\lambda\rfloor(\partial\bar\partial\beta)
=\bar\partial\left(\lambda\rfloor\partial\beta\right)\, ,
\endmatrix$$
showing that also the first two summands are cohomologous to zero.
This proves our contention. \enddemo
\demo{End of the proof of Theorem 1.1}
We prove by contradiction that $V=\emptyset$. In fact, assume that
$x\in V\neq\emptyset$. Then $\Cal I$ contains a nonzero $f$ having
a zero of minimal order $\mu>0$ at $x$. This means that, for holomorphic
coordinates $\zeta_1,\hdots, \zeta_N$ centered at $x$, we have an
expansion
$$f=\dsize\sum_{h=\mu}^\infty{f_h(\zeta)}$$
of $f$ as a convergent series of homogeneous polynomials $f_h$ of
degree $h$, with $f_\mu\neq 0$. Then there is a coordinate $\zeta_j$
such that $\frac{\partial f_\mu}{\partial\zeta_j}\neq 0$. 
By Cartan's Theorem A, since the sheaf of germs of holomorphic vector
fields is coherent, there is a holomorphic vector field $\lambda$ on
$X$ such that $\lambda_x=\left.\frac{\partial}{\partial\zeta_j}\right|_0$.
By Lemma 1.2 we have $\lambda(f)\in\Cal I$. But this gives a contradiction
because $\lambda(f)\neq 0$ has a zero of order $\mu-1$ in $x$.
\par
This completes the proof of the Theorem.
\enddemo
Next we consider the following situation: $F$ is a locally closed subset
of a complex manifold $X$, and $S$ a subset of $F$ which is closed in $F$.
If $Y$ is an open neighborhood of $F$ in $X$ with $\overline{F}\cap Y=F$,
then also $\overline{S}\cap Y=S$. The inclusion map 
$\iota:S\hookrightarrow F$ naturally induces maps:
$\iota^*:H^{p,q}_{\bar\partial}(\frak W_F)@>>>
H^{p,q}_{\bar\partial}(\frak W_S)$,
$\iota_*:H^{p,q}_{\bar\partial}(\frak W_S^{\text{comp}})@>>>
H^{p,q}_{\bar\partial}(\frak W_F^{\text{comp}})$,
$\iota_*:H^{p,q}_{\bar\partial}(\Cal D'_S)@>>>
H^{p,q}_{\bar\partial}(\Cal D'_F)$,
$\iota_*:H^{p,q}_{\bar\partial}(\Cal E'_S)@>>>
H^{p,q}_{\bar\partial}(\Cal E'_F)$
for every $0\leq p,q\leq N$.
We can also consider an open subset $\omega$ of $Y$ and,
corresponding to the inclusion $\sigma:\omega\cap F @>>> F$, the maps
in cohomology:
$\sigma^*:H^{p,q}_{\bar\partial}(\frak W_F)@>>>
H^{p,q}_{\bar\partial}(\frak W_{\omega\cap F})$,
$\sigma_*:H^{p,q}_{\bar\partial}(\frak W_{F\cap\omega}^{\text{comp}})@>>>
H^{p,q}_{\bar\partial}(\frak W_{F}^{\text{comp}}))$,
$\sigma^*:H^{p,q}_{\bar\partial}(\Cal D'_F)@>>>
H^{p,q}_{\bar\partial}(\Cal D'_{\omega\cap F})$,
$\sigma_*:H^{p,q}_{\bar\partial}(\Cal E'_{\omega\cap F})@>>>
H^{p,q}_{\bar\partial}(\Cal E'_F)$
(for $0\leq p,q\leq N$). Denote by $\bold K$ any of the images of the
maps in cohomology considered above. 
Then we obtain, just by repeating the argument of the proof of Theorem 1.1:
\proclaim{Theorem 1.2} With the notation above, if $X$ is a Stein manifold
then $\bold K$ is either zero or infinite dimensional.\endproclaim

\medskip
\noindent
{\bf \S 2 Obstructions to generic embeddings of $CR$ manifolds}
\smallskip
Let $M$ be a smooth (abstract) $CR$ manifold of type $(n,k)$ and let
$\bar\partial_M$ be the tangential Cauchy-Riemann operator on $M$.
Fix a paracompactifying family $\Psi$ of supports in $M$ and consider,
for $0\leq p\leq n+k$, $0\leq q\leq n$, the $\bar\partial_M$-cohomology
groups for smooth differential forms with support in $\Psi$,
denoted by $H^{p,q}_{\bar\partial_M}([\Cal E]^{\Psi}(M))$, and
the corresponding groups for currents with supports in $\Psi$,
denoted by $H^{p,q}_{\bar\partial_M}([\Cal D']^{\Psi}(M))$.
\proclaim{Theorem 2.1} If for some $(p,q)$,
with $0\leq p\leq n+k$, $0\leq q\leq n$,
and a paracompactifying family of supports $\Psi$ in $M$, any one of
the groups $H^{p,q}_{\bar\partial_M}([\Cal E]^{\Psi}(M))$,
$H^{p,q}_{\bar\partial_M}([\Cal D']^{\Psi}(M))$ is finite dimensional
and different from zero,
then there does not exist a {\it generic} $CR$ embedding of $M$ into
any open subset $Y$ of a Stein manifold $X$.
\endproclaim
\demo{Proof}
Assume that $M$ can be generically embedded into an open subset $Y$
of a Stein manifold $X$. The complex dimension of $X$ is $n+k$ and
$M$, being a closed subset of $Y$, is locally closed in $X$.\par
The family $\Phi$ of closed subsets $S$ of $Y$ such that $S\cap M\in\Psi$
is a paracompactifying family in $Y$. \par
In this situation it is a well known consequence of the formal 
Cauchy-Kowalewski theorem (and its dual version) (see [AFN], 
[AHLM], [HN2], [N1], [NV])
that 
$$H^{p,q}_{\bar\partial_M}([\Cal E]^{\Psi}(M))\simeq
H^{p,q}_{\bar\partial}(\frak W_F^{\Phi})\tag \text{$\clubsuit$} $$
and
$$H^{p,q}_{\bar\partial_M}([\Cal D']^{\Psi}(M))\simeq
H^{p,q+k}_{\bar\partial}(\Cal D_M^{\prime\Phi})\, .\tag {$\diamondsuit$} $$
Thus we obtain the conclusion using Theorem 1.1.
\enddemo
\noindent
{\sc Remark} If $M$ has a {\it non-generic} $CR$ embedding 
as a closed submanifold of
an open subset $Y$ of a Stein manifold $X$, then the groups in the
right hand side of ($\clubsuit$) and ($\diamondsuit$) are either
zero or infinite dimensional. But the isomorphism fails, and in fact
the conclusion of Theorem 2.1 is false, as it will be shown by
some examples in the next section.


\medskip
\noindent
{\bf \S 3 Applications}\nopagebreak
\smallskip\nopagebreak
\nopagebreak
1. In particular let $\Omega$ be any domain having a smooth
boundary $M = \partial \Omega$ in an $N$-dimensional Stein manifold $X$. Then
for $0 \leq p \leq N$ and $0 < q \leq N-1$, $H^{p,q}([\Cal E ] (M))$
and $H^{p,q}([\Cal D^{\prime}](M) )$ cannot be finite
dimensional without being zero. They are clearly infinite dimensional for $q =
0$; and if $\Omega \subset\subset X$, we know they are also infinite
dimensional for $q = N-1$ by [HN1].
\smallskip
2. In [Br] it was shown that, if $D$ is a sufficiently small open subset of a
generic $CR$ submanifold $M$ of some open set $\Omega$ in $\Bbb C^{N}$, then
$H^{p,q}([\Cal E ](D))$ is either zero or infinite
dimensional. This follows from Theorem 2.1, without any assumption
on $D$, as far as the embedding $M\hookrightarrow\Omega$ is generic.
\par
Dropping the genericity assumption, the result is still valid for
small open $D$'s because an appropriate holomorphic projection into an
affine $\Bbb C^{n+k}$ will produce a local generic $CR$ embedding.
\smallskip
3. In [Br] it was also pointed out that there exist compact strictly
pseudoconvex $CR$ manifolds $M$ of hypersurface type $(n,1)$, with $n \geq 2$,
which are non-generically $CR$ embedded into some $\Bbb C^{N}$, with $0 < \dim
H^{0,1}([\Cal E ](M)) < \infty$. By Theorem 2, such an $M$
has no generic $CR$ embedding into any Stein manifold. But by [HN1] we know
that the top cohomology groups $H^{p,n}([ \Cal E ](M))$ are
infinite dimensional for $0 \leq p \leq n+1$, due to the fact that $M$ is
embedded, even non-generically, into $\Bbb C^{N}$. Hence there can be no example
of the type pointed out in [Br] with $\dim_{\Bbb R}M = 3$.
\smallskip
4. Suppose $M$ is a compact $CR$ manifold of any type $(n,k)$, $n,k \geq 1$,
which has a non-generic $CR$ embedding in some Stein manifold $X$. Then for $0
\leq p \leq n+k$ the bottom and the top cohomology groups $H^{p,0}( [\Cal E ]
(M))$ and $H^{p,n}([\Cal E ](M))$ are infinite dimensional,
according to [HN1]. Hence finite dimensionality of some bottom or top group
obstructs even non-generic embeddings. In this situation the existence of any
nonzero but finite dimensional intermediate cohomology group $H^{p,q}([\Cal E
](M))$, $0 \leq p \leq n+k$, $0 < q < n$, would obstruct any attempt to make
the non-generic embedding generic. In particular this means that, for such an
$M$, no matter how we embed the Stein manifold $X$ into some $\Bbb C^{N}$, the
$M$ becomes so positioned in $\Bbb C^{N}$ as not to have any one-to-one
holomorphic projection into any affine $\Bbb C^{n+k}$ contained in the $\Bbb
C^{N}$. 
\smallskip
5. Consider a compact smooth orientable $CR$ manifold $M$ of hypersurface
type $(n,1)$, $n \geq 1$, which has a non-generic $CR$ embedding in some
$\Bbb C^{N}$. Then by [HL] there is a holomorphic chain $C$ whose boundary is
$M$ in the sense of currents. Let $V$ denote the support of $C$ and set $F = M
\cup V$. Then $F$ is a closed set in $Y = \Bbb C^{N}$. Hence by Theorem 1, we
have that for $0 \leq p \leq N$ and $0 < q \leq N$, the cohomology groups
$H^{p,q}(\frak W_{F})$ and
$H^{p,q}(\Cal D^{\prime}_{F})$ are either zero or infinite
dimensional.
\par
Suppose $N = n+2$ (so $M$ has real codimension $3$), $n > 1$,
 $M$ is strictly pseudoconvex, and $V$ has only isolated
hypersurface singularities. Then $H_{\bar\partial_M}^{p,q}([\Cal E](M))$
is nonzero and finite dimensional, for $p+q = n$ and 
$0 < q < n$, see [Y]. Thus there are many examples like
the one pointed out by Brinkschulte in [Br]. Moreover this
phenomenon starts to occur as soon as the embedding is in just one complex dimension too high to be generic, so that the embedding would be generic, if the embedding dimension were to be reduced by one.

\bye